\documentclass[12pt,reqno]{article}

\usepackage{amsmath}
\usepackage{amsthm}
\usepackage{amscd}
\usepackage{amssymb}
\usepackage{graphicx}
\usepackage{epsf}
\theoremstyle{plain}
\newtheorem{theorem}{Theorem}[section]

\newtheorem{corollary}[theorem]{Corollary}

\newtheorem{definition}[theorem]{Definition}
\newtheorem{remark}[theorem]{Remark}
\newtheorem{example}[theorem]{Example}

\def\less{<}

\def\int{\text{int}}
\def\invlimit{\smash{\lim\limits_{\raise1pt\hbox{$\longleftarrow$}}}\vphantom{\big(}}
\def\inter{\hskip 1.5pt\raise4pt\hbox{$^\circ$}\kern -1.6ex}
\def\skel(#1,#2){#1^{(#2)}}
\def\hyp {\hbox {\rm {H \kern -2.8ex I}\kern 1.25ex}}
\def\reals {\hbox {\rm {R \kern -2.8ex I}\kern 1.15ex}}
\def\integers {\hbox {\rm { Z \kern -2.8ex Z}\kern 1.15ex}}
\def\naturals {\hbox {\rm {N \kern -2.8ex I}\kern 1.20ex}}
\def\rationals {\hbox {\rm { Q \kern -2.2ex l}\kern 1.15ex}}
\def\hyp {\hbox {\rm {H \kern -2.7ex I}\kern 1.25ex}}

\begin{document}
\title{On boundary primitive manifolds and a theorem of Casson-Gordon}

\author{ Yoav Moriah\thanks{Supported by The Fund for Promoting Research at the
Technion,  grant 100-127 and the Technion VRP fund, grant 100-127.}}

\date{}

\maketitle

\begin{abstract} In this paper it is shown that manifolds admitting minimal genus weakly reducible
but irreducible Heegaard splittings contain an essential surface. This is an extension of a well known 
theorem of Casson-Gordon to manifolds with non-empty boundary. The situation for non-minimal 
genus Heegaard splittings is also investigated and it is shown that boundary stabilizations are
stabilizations for manifolds which are boundary primitive.

\quote{{\it Keywords } Heegaard splittings, weakly reducible, irreducible, primitive meridian

{\it 2000 AMS classification } 57N25 }

\end{abstract}
\section{Introduction}

\label{intro}

\vskip10pt

A well known result of Casson-Gordon  (see [CG] Theorem 3.1) states that if $M$ is a \underline {closed}
orientable manifold and $(V_1,V_2)$ an irreducible but weakly reducible Heegaard splitting of $M$ then
$M$  contains an essential surface of positive genus.  This theorem is extremely useful and it was a natural 
thing to expect an extension of it to manifolds with boundary. Surprisingly, so far, the statement of [CG]  
Theorem 3.1 does not extend as is to manifolds with boundary and the emerging picture is rather
complicated as will become clear from the following theorems. First the positive results:

\vskip10pt

\noindent {\bf Theorem 3.1:} {\it Let $M^3$ be an orientable $3$-manifold which has a weakly reducible 
Heegaard splitting of minimal genus, then $M$ contains an essential surface of positive genus.}

\vskip10pt

This raises the question of what can be said about irreducible but weakly reducible Heegaard splittings which
are not of minimal genus. We have:

\vskip10pt 

\noindent {\bf Theorem \ref{mingengcor}:} {\it Let $M^3$ be an orientable $3$-manifold with a single 
boundarycomponent of genus $h$. Assume that $M$ has a weakly reducible but irreducible Heegaard  
splitting $(V_1, V_2)$ of genus $g$, then either $M$ contains an essential surface or $(V_1, V_2)$ is 
a boundary stabilization and $M$ has a strongly irreducible Heegaard splitting of genus $g - h$.}

\vskip10pt
 
An extension of the Casson-Gordon theorem to non-minimal genus Heegaard splittings of manifolds with
boundary can fail in only one way i.e., if there is a manifold with an irreducible but weakly reducible
non-minimal Heegaard splitting, which contains only boundary parallel incompressible surfaces. Such a surface
will separate a region homeomorphic to $\mbox {(boundary)} \times I$ and the Heegaard splitting will induce a
{\it  standard} Heegaard splitting on this region (otherwise the original Heegaard splitting will be reducible by
[ST]). 

A somewhat simple example for the failure of the Casson-Gordon theorem to manifolds with boundary
is given in Example 6.1 of [Se], where a genus three weakly reducible and irreducible Heegaard splitting for 
the  complement of the three component trivial open chain link is presented. Since the complement of the link is
homeomorphic to a $\mbox{(pair of pants)} \times S^1$ it contains no closed essential surfaces. The 
Heegaard splitting in the example has all three boundary components contained in one compression body.
However a minimal genus splitting for $\mbox{(pair of pants)} \times S^1$ is of genus two with one 
compression body containing two boundary components and the other one boundary component.
The question is still open for manifolds with two or less boundary components.

We can take the opposite point of view: start with an irreducible Heegaard splitting and add a $(boundary
\times I)$ with a  standard Heegaard splitting to it. This operation will be called a {\it boundary stabilization}.
This does not change the manifold but will give a new amalgamated  Heegaard splitting which if weakly
reducible will give rise to an incompressible but boundary parallel surface and thus a candidate for a counter
example to an extension of the Casson-Gordon theorem. In order to do this the notion of a  
{\it $\gamma$-primitive} Heegaard splitting is defined in Section 4. This approach also runs into difficulty as we
have:

\vskip10pt

\noindent {\bf Theorem \ref{prim}:} {\it If an orientable $3$-manifold $M$ has a $\gamma$-primitive 
Heegaard splitting then every boundary stabilization on the component containing the curve $\gamma$ is a
stabilization. }

\vskip6pt

Since  $\gamma$-primitive Heegaard splittings are very common weakly reducible but irreducible 
non-minimal genus Heegaard splittings are hard to find and hence possible counter examples to 
the remaining cases are also hard to find.

\vskip10pt

\noindent{\bf Remark:} In [LM] Theorem 1.3 we stated an extension of this result to manifolds with boundary
but  unfortunately  the statement and the proof given there are not quite right (see the footnote in Section 2). I
would like to thank  T.  Kobayashi for  pointing this out to me. To the best of my knowledge no such extensions
of [CG] Theorem 3.1 appeared before [LM], and the problems mentioned above do not affect the other results of
that paper.  In this paper the situation is corrected. 

Extensions to [CG] Theorem 3.1 have been proved by E. Sedgwick [Sc] who proves a similar theorem 
to  Theorem 3.1 and recently by  T. Kobayashi in [Ko] Proposition 4.2 which states that, if $M$ 
has a weakly reducible Heegaard splitting  then either the Heegaard splitting is reducible or $M$ contains
incompressible surfaces of positive genus (which might be boundary parallel i.e., not essential). The theorems
presented above are  a strengthening of this result.

\vskip10pt

\section{Preliminaries}

\vskip10pt

\label{prel}
In this paper it is assumed that all manifolds and surfaces will be orientable unless otherwise specified.

A compression body  $V$  is a  compact orientable and connected 3-manifold with a preferred 
boundary component $\partial_+V$  and is obtained from a collar of $\partial_+ V$ by attaching 
2-handles and 3-handles, so that the connected components of  $\partial_- V$ = $\partial V - \partial_+ V$ 
are all distinct from  $S^2$.  The extreme cases, where  $V$  is a handlebody i.e., $\partial_- V = \emptyset$,
or where $V = \partial_+V \times I$, are allowed.  Alternatively we can think of $V$ as obtained from
$(\partial_-V) \times I$ by attaching $1$-handles to $(\partial_-V) \times \{1\}$. An annulus in a 
compression body will be called a {\it vertical (or a spanning) annulus } if it has its boundary components
on different boundary components of the compression body.

Given a manifold  $M^3$ a {\it Heegaard splitting } for $M$ is a decomposition  $M = V_1 \cup V_2$ into
two compression bodies $(V_1, V_2)$  so that $V_1 \cap V_2 = \partial V_1 = \partial V_2 = \Sigma$ The
surface $\Sigma$ will be call the {\it Heegaard splitting surface}. 

A  Heegaard splitting $(V_1,V_2)$ for a manifold $M$ will be called {\it reducible } if there are
essential disks$D_1 \subset V_1$ and $D_2 \subset V_2$ so that $\partial D_1 = \partial D_2 \subset \Sigma$.

A  Heegaard splitting $(V_1,V_2)$ for a manifold $M$ will be called {\it weakly reducible } if there are  disjoint
essential disks$D_1 \subset V_1$ and $D_2 \subset V_2$. Otherwise it will be called {\it strongly irreducible}.

Let $M$ be a $3$-manifold which is homeomorphic to a $\mbox{(surface)} \times I $. A Heegaard 
splitting $(V_1,V_2)$ of $M$ will be called {\it standard} if it is homeomorphic to one of the following types:
\begin{itemize}
\item[(I)] $V_1 \cong (surface) \times [0, \frac1 2]$, $V_2 \cong (surface) \times [ \frac1 2, 1]$ and
$ \partial_+ V_1 =  \partial_+ V_2 = (surface) \times [ \frac1 2] $
\item[(II)] If $ \{p\} \in (surface) $ is a point then for $0 \less\epsilon \less \frac1 2$
 
\obeylines 
$V_1 \cong ((surface) \times [0, \epsilon]) \cup( N(p)\times I) \cup  ((surface) \times [1- \epsilon,1])$ and 
$V_2 = cl(M - V_1)$. 
\end {itemize}

Note that $V_2$ is a regular neighborhood  of a once punctured surface and hence is a handlebody and
$V_1$ is a compression body with one boundary component $\partial _+$ of genus $2g$ and two   
boundary components $\partial _-$ of genus $g$, where $g = genus (surface)$. In [ST] it is proved that any
irreducible Heegaard splitting of  $(surface) \times I $ is homeomorphic to one of the above two types.

A closed surface $F \subset M$ will be called {\it essential} if it incompressible and non-boundary parallel.

Given a closed (possibly disconnected) surface $\Sigma \subset M$ and a system of  pairwise disjoint 
non-parallel compressing disks $\Delta$ for $\Sigma$ define (as in [LM]) 
$\Sigma_0 = \sigma(\Sigma,\Delta)$ to be the surface obtained from $\Sigma$ by compressing along
$\Delta$. Let $c(\Sigma) = \sum_i (1 - \chi(\Sigma_i))$, where the sum is taken over all components
$\Sigma_i$ of $\Sigma$ which are not $2$-spheres. The complexity of the system $\Delta$ is defined
to be:
\vskip8pt
\hskip30pt $c(\Delta) = c(\Sigma) - c(\Sigma_0)$

\vskip8pt

For a given Heegaard splitting surface $\Sigma$ for $M$  we will assume that  a system of compressing
disks $\Delta = \Delta_1 \cup \Delta_2$, where $\Delta_i \subset V_i  $,  satisfies:

\begin{itemize}
\item[(a)] $\Delta_i \neq \emptyset$ for both $i = 1,2$.  i.e., $\Delta$ contains disks on both sides of $\Sigma$.
\item[(b)] $\Delta$ is maximal with respect to $c(\Delta)$ over all systems $\Delta$ satisfying (a).
\end {itemize}

\begin{definition}
\label{induced}\rm

Let $\Sigma^*$ be the surface $\Sigma_0$ less the $2$-sphere components and the components which are
contained in $V_1$ or $V_2$  \footnote {This definition of $\Sigma^*$ is different from that of [LM] in that 
it excludes the components of $\Sigma_0$ which are contained in $V_1$ or $V_2$. The problems 
in the proof of Theorem 1.3 of [LM] emanate from the fact that with the definition given there the 
component N is not correctly defined. In particular, note that with this modified definition the mistake 
in the proof of  Theorem 1.3  of  [LM]  disappears. However, as expected, this  will not correct the 
given proof: The point is that, with the modified definition, the statement of Lemma 1.2 (b) becomes
wrong  (compare also [Ko]).}. 
Let $N_0$ denote the  closure of a component of $M - \Sigma_0$ which is not a $3$-ball and let $N$ 
denote the closure of a component of   $M - \Sigma^*$ which contains $N_0$.  By the symmetry 
between $V_1$  and $V_2$ we can assume that $N_0 \subset V_1 \cup \eta(\Delta_2)$. Now set 
$U_1 = (V_1 \cap N_0) - \eta(\Sigma \cup \Delta)$ and $U_2 = N - U_1$. By Lemma 1.2 (a) of 
[LM] the pair $(U_1,U_2)$ is a Heegaard splitting for $N$ and will be called {\it  the induced 
Heegaard splitting on $N$}.

\end{definition}

 If $(V_1,V_2)$ is an irreducible Heegaard splitting of $M$ then
$(U_1,U_2)$ is an irreducible Heegaard splitting of $N$ by Lemma 1.2(c) of [LM].

 Given two manifolds $M_1$ and $M_2$ with respective Heegaard splittings $(U^1_1,U^1_2)$ and 
$(U^2_1,U^2_2)$, assume further that there are homeomorphic boundary components 
$F_1 \subset \partial\_ U^1_1$ and  $F_2 \subset \partial\_ U^1_2$.  Let $M$ be a manifold obtained by
gluing $F_1$ and $F_2$.  We can obtain a Heegaard splitting $(V_1,V_2)$ for $M$ by a process called
{\it amalgamation} (see [MS] and [Sc]). The process of amalgamation reconstructs the original Heegaard splitting 
$(V_1,V_2)$ of $M$ from the Heegaard splittings induced on the components $N_i$ of $M - \Sigma^*$
(see Fig.1).


\vbox{{\epsfysize110pt\epsfbox{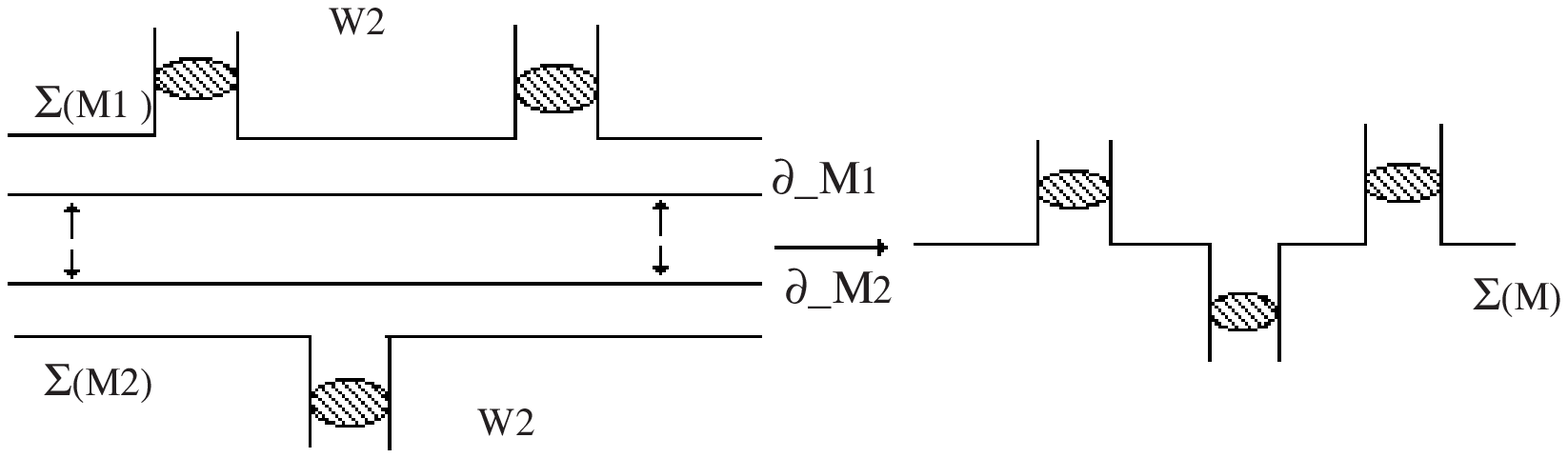}}

\vskip40pt

\centerline{Fig.1}}

\vskip35pt

\section{Essential surfaces}

\vskip25pt

In this section we prove two of the main theorems:

\vskip10pt

\begin{theorem} 
\label{mingen}
Let $M^3$ be an orientable $3$-manifold which has a weakly reducible Heegaard 
splitting of minimal genus, then $M$ contains an essential surface of positive genus.

\end {theorem}

\vskip10pt

\begin{proof} 
Let $(V_1, V_2)$ be an irreducible and weakly reducible Heegaard splitting of minimal genus for $M$ and
$\Sigma, \Delta, \Sigma_0$ and  $\Sigma^*$ be as above. Since $\Sigma$ is connected there must be at 
least one component $S$ of $\Sigma_0$ so that both of $S \cap  V_1$ and $S \cap V_2$ are not empty. 
Since the Heegaard splitting is minimal it is irreducible and since the surface $S$ contains disks from both
$V_1$ and 
$V_2$ it is not a $2$-sphere and hence is in $\Sigma ^*$. By the proof of Theorem 3.1  in [CG] since $\Delta$ is 
maximal then $S$ is incompressible. 

It remains to show that $S$ is not boundary parallel when $\partial M \neq \emptyset$. If $S$ is boundary 
parallel  then $M - S$  has two  components and the closure of one $M_1$, is homeomorphic to $S \times I$. 
Let $M_0$ be the  closure of  the other component. Note that  $M_0$ is homeomorphic to $M$ and that  
$(V_1, V_2)$ induces a Heegaard splitting on both of $M_0$ and $M_1$ as in Definition \ref{induced}. 
Assume that the component of 
$\partial M$ homeomorphic to $S$ is contained in, say, $V_1$. Since $S \cap V_1 \neq \emptyset$ then  
$S \cap V_1$ must consist of a single disk as otherwise the induced Heegaard splitting on $M_1$ will be 
reducible: As all Heegaard splittings of $S \times I$ are standard and since $S \cap V_1 \neq \emptyset$ then the
induced Heegaard splitting is of type II in the terminology of [ST] and is reducible if there is more than one such
disk. Furthermore if the induced Heegaard splitting on $S \times I$ is reducible then by [LM] Lemma 1.2(c) ( see
also [Ko] Lemma 4.6)  it follows that $(V_1, V_2)$ is reducible in contradiction.

This implies that the genus of the induced Heegaard splitting on $M_1$ is $2 \times g(S)$. The formula for
computing the genus of the amalgamated  Heegaard splitting from the genus of the Heegaard splitting of the
components $M_0, M_1$ is ${ \bar g}(M) = { \bar g}(M_1) + { \bar g}(M_0) - g(S)$ 
where ${ \bar g}$ is the genus of the induced Heegaard splitting (see [Sc]). Hence ${ \bar g}(M_0)$
must be strictly smaller than ${ \bar g}(M)$ in contradiction to the fact that $M_0$ is homeomorphic to $M$
and that $(V_1, V_2)$   is a minimal genus Heegaard splitting  of $M$.

\end{proof}

\begin{theorem} 
\label{mingeng}
Let $M^3$ be an orientable $3$-manifold with a single boundary component of genus $h$. Assume that $M$
has a weakly reducible but irreducible Heegaard  splitting $(V_1, V_2)$ of genus $g$, then either $M$ contains
an essential surface or $M$ has a strongly irreducible Heegaard splitting of genus $g - h$.

\end {theorem}
\vskip10pt
\begin{proof}  Let $(V_1, V_2)$ be an irreducible and weakly reducible Heegaard splitting of  genus $g$
for $M$ and $\Sigma, \Delta, \Sigma_0, \Sigma^*$ and  $S$ be as above.  If $S$ is essential we are done . 
So we assume that $S$ is boundary parallel and hence $g(S) = h$. As in the proof of Theorem \ref{mingen}
let $M_0$ and $M_1$ be the closure of the  components of $M - S$. The Heegaard splitting $(V_1, V_2)$
of $M$ induces a Heegaard splitting of genus $2h$ on $M_1 = S \times I$ and a Heegaard splitting 
$(U_1, U_2)$ of genus $g - h$ on $M_0$ as in Definition \ref{induced}.

The Heegaard splitting $(U_1,U_2)$ is irreducible as otherwise  it follows from [LM] Lemma 1.2(c) 
(see also [Ko]  Lemma 4.6)  that $(V_1, V_2)$ is reducible in contradiction. If $(U_1, U_2)$ is
weakly reducible then by Theorem \ref{mingen} $M_0$ has an incompressible surface $S_0$. If
$S_0$ is essential in $M_0$ we are done since $M \cong M_0$. If $S_0$ is boundary parallel then
since $M_0$ has a single boundary component which is homeomorphic to $S$ we have two
incompressible surfaces $S$ and $S_0$ which are boundary parallel. Hence the closure of the
component of $M - S_0$ is homeomorphic to 
$S_0 \times I \cong S_0 \times [0,{\frac 1 2} ]\cup_S S \times [{\frac 1 2},1 ]$. The amalgamation of
the genus $2h$ Heegaard splittings of $S_0 \times [0,{\frac 1 2}]$ and $S \times [{\frac 1 2},1 ]$ will
induce a Heegaard splitting of genus $2h + 2h - h$ on $S_0 \times I$ which is reducible as all Heegaard
splittings of a $(surface) \times I$ are standard by [ST]. But this implies as before that $(V_1, V_2)$ is 
reducible in contradiction. Hence $(U_1, U_2)$ is a strongly irreducible Heegaard splitting of genus
$g - h$ of $M \cong M_0$ (see Fig.2).

\end{proof}


\vbox{{\epsfysize95pt\epsfbox{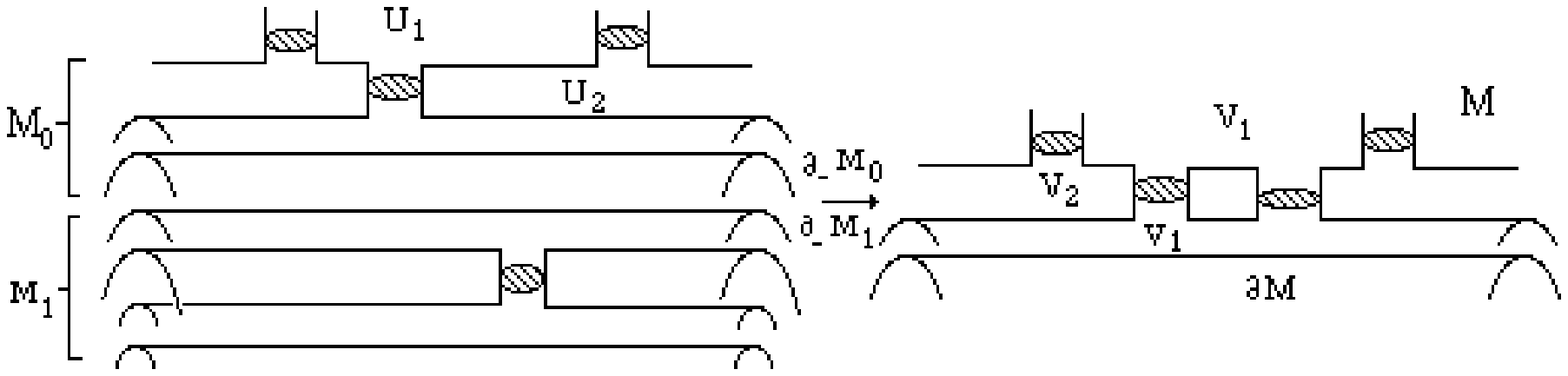}}

\vskip15pt \centerline {Fig.2}}

\vskip30pt

\begin{remark}\rm The same argument as in the above proof could be used for manifolds with more than
one boundary component. However the statement of the theorem in that case would be very
complicated.
\end{remark}

\vskip20pt

\section{Boundary primitive manifolds}

\vskip 20pt

Given a manifold $M$ with boundary components $\partial M^1 ,\dots, \partial M^k$ of corresponding 
genus $g^1 , \dots , g^k$ and a Heegaard splitting $(V_1,V_2)$ for $M$ of genus $g$ we can always 
obtain a new Heegaard splitting $(U_1^i,U^i_2)$ of genus $g +g^i,  i = 1, \dots , k$, by gluing a 
$\partial M^i \times I$  to  the $\partial M^i$ boundary component and then amalgamating the standard 
Heegaard splitting of genus $2g^i $ of $\partial M^i \times I$ with the given Heegaard splitting 
$(V_1,V_2)$ of $M$ (as indicated  in Fig.1). 

\begin{definition}
\label{stabdef}
\rm  The construction above will be called {\it boundary stabilization on the i-th
boundary component}. If there is a single boundary component or no ambiguity we can just
use {\it boundary stabilization}.
\end{definition}

\vskip10pt

\noindent We can now restate a stronger form of Theorem \ref{mingeng} namely:

\vskip10pt

\begin{theorem} 

\label{mingengcor}
Let $M^3$ be an orientable $3$-manifold with a single boundary component of genus $h$. Assume that $M$
has a weakly reducible but irreducible Heegaard  splitting $(V_1, V_2)$ of genus $g$, then either $M$ contains 
an essential surface or $(V_1, V_2)$ is a boundary stabilization and $M$ has a strongly irreducible Heegaard 
splitting of genus $g - h$.

\end {theorem}
\qed

\vskip10pt

\begin{remark}\rm
The Heegaard splitting $(U_1^i,U^i_2)$ is clearly weakly reducible by the construction and the 
question arises of  when is it irreducible?  This question is of interest as it was shown in [LM1] 
that it is relatively easy to find manifolds with an arbitrarily large number of strongly irreducible 
Heegaard splittings. It is much more difficult to find manifolds with irreducible but weakly reducible
Heegaard splittings. 
\end{remark}

\vskip6pt

We say that an element $x $  in a free group $F_n$ is {\it primitive} if it belongs to some basis for 
$F_n$. A curve on a handlebody $H$  is {\it primitive} if it represents a primitive element in the free
group $\pi_1(H)$. An annulus $A$ on  $H$ is {\it primitive} if its core curve is primitive. Note that
a curve on a handlebody is primitive if and only if there is an essential disk in the handlebody 
intersecting the curve in a single point.

\vskip8pt

\begin{definition}{\rm Let $M$ be a $3$-manifold with incompressible boundary components
$\partial M^1 ,\dots, \partial M^k$ .  Let $\gamma \subset \partial M^i$  be an essential simple 
closed curve. A Heegaard splitting $(V_1,V_2)$ of  $M$  will be called {\it  $(\gamma-primitive)$} 
if there is an  annulus $A$ in $V_1$ or $V_2$, say $V_1$, with $\gamma$ as one boundary 
component of $A$ and  the other a curve on the Heegaard surface $\Sigma$ which 
intersects  an essential disk of $V_2$ in a single point}.
\end {definition}

\vskip8pt

\begin{theorem}
\label{prim}
If an orientable $3$-manifold $M$ has a $\gamma$-primitive Heegaard splitting then every boundary
stabilization on the component containing the curve $\gamma$ is a stabilization. 

\end{theorem}

\vskip8pt

\begin{proof} Let $\partial M^i$ be the boundary component on which we are going to stabilize.
Assume that $genus (\partial M^i) = g_i$, hence we amalgamate the given  genus Heegaard splitting
$(W_1, W_2)$ of $M$ with a genus $2g_i$ Heegaard splitting $(U_1, U_2)$ of $\partial M^i \times I$. 
The Heegaard splitting $(U_1, U_2)$ is standard of type II in the terminology of [ST].

Since the Heegaard splitting of $M$ is $\gamma$-primitive there is some curve $\gamma \subset \partial M^i$
which  co bounds an annulus $A'$ so that the other boundary component of $A'$ meets an essential disk  $D_2$ of
$W_2$, say, in a single point. When $\partial M^i$ is identified with  $\partial M^i \times I$ the curve $\gamma$
determines an annulus $ A = \gamma \times I \subset \partial M^i \times I$.  Assume that the handlebody
component of the standard Heegaard $(U_1, U_2)$ splitting of  $\partial M^i \times I$ is $U_2$. After
an ambient isotopy of $\partial M^i \times I$ we can always assume that the vertical arc $\{p \} \times I$ is
contained in $A$ hence $A \cap U_2$ is an essential disk $D_1$ in $U_2$.   In the process of amalgamating 
the Heegaard splittings the handlebody $U_2$ gets glued to $W_1$ and the compression body $U_1$ gets
glued to $W_2$. It is possible that $D_1$  and $D_2$ will get identified with disks which are not properly
embedded. However this can be corrected by a small isotopy. Now the two disks $D_1 \subset W_1$ and $D_2
\subset W_2$ still intersect in a  single point. In the amalgamated Heegaard splitting $(V_1,V_2)$ of $M$ 
we have that $D_1 \subset V_1$ and $ D_2 \subset V_2$ and hence it is a  stabilization (see
Fig. 3). 

\end{proof}

\begin{definition}{\rm Let $M$ be a $3$-manifold with incompressible boundary components 
$\partial M^1, ..., \partial M^k$. The manifold $M$ will be called {\it boundary primitive or
$(\partial-primitive)$}  if for each Heegaard splitting of minimal genus of $M$ and for each boundary
component $\partial M^i$ there is some curve $\gamma \subset \partial M^i$ for which the Heegaard splitting is
$\gamma-primitive$.  

In particular if $M = S^3 - N(K)$ where $K \subset S^3$ is a knot and $\gamma = \mu$ is
a meridian  curve  we will say that the  Heegaard splitting is $\mu$-primitive and that $E(K)$ is $\mu$-primitive
if all its Heegaard splittings of minimal genus are  $\mu$-primitive}. 

\end{definition}

\vskip8pt

\begin{corollary}
\label{primcor}
If an orientable $3$-manifold $M$ is $\partial$-primitive then every boundary stabilization of a minimal
genus Heegaard splitting is a stabilization. 
\end{corollary}
\qed

\begin{corollary} Let $K \subset S^3$ be  a $\mu$-primitive knot. Assume that $g(E(K)) = g$ and
that $E(K)$ has an irreducible but weakly reducible Heegaard splitting of genus $g + 1$ then $E(K)$
contains an essential surface.
\end{corollary}
\qed


\vbox{{\epsfysize120pt\epsfbox{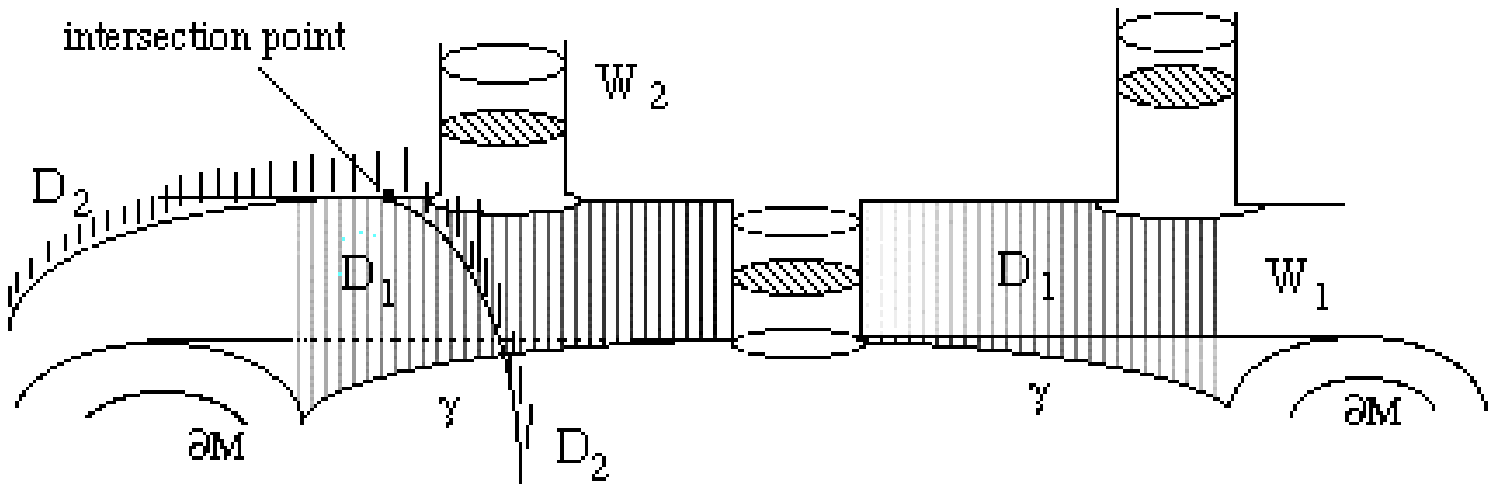}}
\vskip25pt\centerline{Fig.3}}

\vskip10pt
\vskip15pt

\begin{remark}\rm Knot complements $E(K_i), i = 1,2$, in $S^3$ which are tunnel number super additive, 
i.e., $t(K_1 \# K_2) = t(K_1) + t(K_2) +1$ are examples of manifolds which are not $\mu$-primitive. 
These knots exist by independent results of Moriah, Rubinstein [MR] and Morimoto, Sakuma and Yokota [MSY].

\end{remark}

\begin{example} \rm In [Ko] it is proved that all $2$-bridge knots are $\mu$-primitive. 
\end{example}

\vskip10pt

\begin{example} \rm Knots in $S^3$ which admit a $(g,1)$-decomposition (see e.g. [MS1]) have a 
Heegaard splitting which is $\mu$-primitive.
\end{example}

\vskip10pt

\begin{example} \rm Any Heegaard splitting of a knot in a $2n$-plat projection which is induced
by a top or bottom minimal tunnels  is $\mu$-primitive (see [LM]).
\end{example}

\vskip10pt

\begin{remark} \rm All knots $K \subset S^3$ have $g(E(K) +1$ $\mu$-primitive Heegaard splittings 
obtained from a minimal genus Heegaard splitting by stabilizing once  i.e., by drilling a small tunnel from 
$V_2$ and adding it  to $V_1$ so that the new $1$-handle of $V_2$ intersects a vertical annulus of 
$V_1$ in a single point.

\end{remark}

\vskip10pt

\begin{remark}\rm Note that every Heegaard splitting which is a boundary stabilization is $\mu$-primitive.
As the disk $D_1$ in $V_1$ will intersect a curve bounding a vertical annulus with a meridian in $V_2$
in a single point.
\end{remark}
\vskip10pt

\begin{remark}\rm
The previous definitions and examples raise the question of whether a given knot complement can have
Heegaard splittings which are $\mu$-primitive and others which are not. The answer to this question 
is affirmative and examples of such knots are torus knots. Given a torus knot $K(p,q) \subset S^3$, 
where $g.c.d.(p,q) = 1$, then $E(K(p,q))$ has three genus two Heegaard splitting  if and only if 
$p \neq \pm1 \hskip2pt mod  \hskip2pt q$ and $q \neq \pm1 \hskip2pt mod  \hskip2pt p$  (see [Mo],
[BRZ]). One of the Heegaard splittings is obtained by considering a decomposition of $E(K(p,q))$ into 
two solid tori glued along a
$p,q$ - annulus and then drilling out a neighborhood of an essential arc from this annulus. Thus the
two tori are glued along a disk and form a genus two handlebody the spine of which is composed of
two loops $x$ and $y$.  The complement of this  handlebody is homeomorphic to 
$T^2 \times I \cup N(essential \hskip3pt arc)$ which is a genus two compression body. The loops 
$x$ and $y$ generate the following presentation  for the fundamental group of $E(K(p,q))$: 
$\pi_1(E(K(p,q)) \cong  \hskip3pt < x, y  \hskip2pt|  \hskip2pt x^p = y^q >$. Since  $g.c.d.(p,q) = 1$ 
we can find positive integers  $r < q, s < p$ so that $rp - sq =1$. A curve $\mu$ representing the element 
$x^s y^r$  is a meridian of $E(K(p,q))$ (see for example  Proposition 3.28 of [BZ]). Now choose $p$ 
and $q$ so that $min\{r, s\} \geq 2$, The main theorem of [CMZ] shows that $\mu$  is 
not primitive in $F(x,y)$ and in particular this Heegaard splitting is not $\mu$-primitive. 
However the other two Heegaard splittings are $\mu$-primitive as shown in  
[Mo] and [BRZ] also in [MS1].

\end{remark}

\vskip25pt

\section {References}

\vskip20pt

\noindent[BZ] \hskip28pt G. Burde, H. Zieschang; {\it knots}, de Gruyter  Studies in Math.  5 

\noindent  \hskip54pt Berlin, New York 1985.

\vskip10pt

\noindent[BRZ] \hskip20pt M. Boileau, M. Rost H. Zieschang; {\it On Heegaard decompositions 

\noindent  \hskip54pt  of torus knots exteriors and related Seifert fibered Spaces}, 

\noindent  \hskip54pt Math. Ann. 279 (1988) 553 - 581.

\vskip10pt

\noindent[CG] \hskip28pt A.Casson, C. Gordon; {\it Reducing Heegaard splittings}, Topology and 

\noindent  \hskip54pt  its Applications 27, (1987) 275 - 283 .

\vskip10pt

\noindent [CMZ] \hskip18pt M. Cohen, W. Metzler, B. Zimmermann; {\it What does a  basis of

\noindent  \hskip54pt   F(a,b) look like?}, Math. Ann .257, (1981) 435 - 445 .

\vskip10pt 

\noindent [Ko] \hskip28pt T. Kobayashi; {\it Heegaard splittings of exteriors of two bridge knots}, 

\noindent  \hskip54pt preprint. 

\vskip10pt

\noindent [LM] \hskip28pt M. Lustig, Y. Moriah; {\it Closed incompressible surfaces in comple-

\noindent  \hskip54pt ments of wide knots and links}, Topology and its Applications 

\noindent  \hskip54pt 92 (1999) 1 - 13.

\vskip10pt

\noindent [LM1] \hskip22pt M. Lustig, Y. Moriah; {\it 3-manifolds with irreducible Heegaard 

\noindent  \hskip54pt splittings of high genus}, 		Topology,  39  (2000),  589 - 618 .

\vskip10pt

\noindent [Mo] \hskip28pt  Y. Moriah; {\it Heegaard splittings of Seifert fibered spaces}, Invent.

\noindent  \hskip56pt Math. 91 (1988),  465 - 481 .

\vskip10pt

\noindent [MR] \hskip28pt  Y. Moriah, H. Rubinstein; {\it  Heegaard structure of negatively curved 

 \hskip42pt 3-manifolds},      Comm. in  Ann. and Geom. 5 (1997),  375 - 412 .

\vskip10pt

\noindent [MS] \hskip28pt Y. Moriah, J. Schultens; {\it Irreducible Heegaard splittings of Seifert 

\noindent \hskip54pt fibered spaces are either vertical or horizontal}, Topology 37, (1998)  

\noindent \hskip54pt 1089 - 1112 .

\vskip10pt

\noindent [MS1] \hskip20pt K. Morimoto, M. Sakuma; {\it On unknotting tunnels for knots},

\noindent \hskip54pt Math. Ann. 289, (1991) 143 - 167.

\vskip10pt

\noindent [MSY] \hskip20pt K. Morimoto, M. Sakuma, Y. Yokota; {\it Examples of tunnel number   

\hskip36pt  one knots which have the property  "1 + 1 = 3"} Math. Proc. 

\hskip36pt Cambridge Phil. Soc. 199, (1996)   113 - 118.

\vskip10pt

\noindent [Sc] \hskip30pt J. Schultens; {\it The classification of  Heegaard splittings for  

\noindent \hskip54pt (compact orientable surface)$\times S^1$}, Proc. London Math. Soc. 67,   

\noindent \hskip54pt (1993) 425 - 448 .

\vskip10pt

\noindent [Se] \hskip30pt E. Sedgwick; {\it  Genus two 3-manifolds are built from handle number 

\noindent \hskip54pt one pieces},  preprint.

\vskip10pt

\noindent [ST] \hskip28pt M. Scharlemann, A. Thompson; {\it  Heegaard splittings of 

\noindent \hskip54pt $(surfaces) \times I$ are standard}, Math. Ann. 95 (1993), 549 - 564.

\vskip65pt

\obeyspaces   Yoav Moriah

\obeyspaces   Department of Mathematics

 \obeyspaces  Technion, Haifa  32000,

\obeyspaces  Israel

 ymoriah@tx.technion.ac.il

\end{document}